\documentclass{amsart}
\usepackage{amssymb,latexsym}
\theoremstyle{plain}
\newtheorem{theorem}{Theorem}
\newtheorem{corollary}{Corollary}

\newtheorem{lemma}{Lemma}
\theoremstyle{definition}
\newtheorem{definition}{Definition}

\begin{document}

\title{Koszul cohomology and singular curves}
\author{E. Ballico, C. Fontanari, L. Tasin}
\address{Department of Mathematics\\
University of Trento\\
Via Sommarive 14\\
38123 Povo (TN), Italy}
\email{\{ballico,fontanar\}@science.unitn.it, luca.tasin.1@studenti.unitn.it}
\thanks{The authors are partially supported by MIUR and GNSAGA of INdAM (Italy).}

\begin{abstract}
We investigate Koszul cohomology on irreducible nodal curves following 
the lines of \cite{a}. In particular, we prove both Green and 
Green-Lazarsfeld conjectures for any $k$-gonal nodal curve 
which is general in the sense of \cite{bf1}.
\end{abstract}

\maketitle

\section{Introduction}

Let $X$ be a complex projective curve. For any line bundle 
$L\in \mbox{Pic}(X)$ and all integers $p, q$, let $K_{p,q}(X,L)$
denote the Koszul cohomology groups introduced in \cite{g} as 
the cohomology of the complex:
$$
\wedge^{p+1}H^0(L) \otimes H^0(L^{q-1}) \to \wedge^p H^0(L) 
\otimes H^0(L^{q}) \to \wedge^{p-1}H^0(L) \otimes H^0(L^{q+1}).
$$

Green's conjecture states that $K_{p,1}(X,\omega_X)=0$ if and only 
if $p \ge g-\mathrm{Cliff}(X)-1$, where $\mathrm{Cliff}(X)$ 
is the Clifford index of $X$, while Green-Lazarsfeld conjecture 
(see \cite{gl}, Conjecture (3.7)) predicts that for every line bundle 
$L$ on $X$ of sufficiently large degree $K_{p,1}(X,L)=0$ if and only if 
$p \ge r-\mathrm{gon}(X)+1$, where $r$ is the (projective) dimension of $L$ 
and $\mathrm{gon}(X)$ is the gonality of $X$. 

Both Green and Green-Lazarsfeld conjectures have been verified 
for the general curve of genus $g$ 
(see \cite{v1}, \cite{v2}, \cite{av}, \cite{ap}) and for the general 
$d$-gonal curve of genus $g$ (see see \cite{t} for $d \le g/3$, \cite{v1}, 
Corollary 1 on p. 365, for $d \ge g/3$, \cite{av}, \cite{a}). 

Indeed, \cite{a} shows that Green's conjecture is satisfied 
for any smooth $d$-gonal curve verifying a suitable linear 
growth condition on the dimension of Brill-Noether varieties 
of pencils which holds for the general $d$-gonal curve.
The arguments in \cite{a}, taking the path opened in \cite{v1}, 
rely on suitable degenerations to irreducible nodal curves. 
As a by-product, they imply that a general irreducible nodal 
curve $Y$ of genus $g=2k+1$ has not extra-syzygies, i.e. 
$K_{k,1}(Y, \omega_Y)=0$ (see \cite{a}, Proposition 7).   

Here instead we regard singular curves not just as a powerful tool 
but as a natural geometric object being interesting also in its own 
and we push further the intuition underlying \cite{a}, proof of Theorem 2
(see also \cite{bf2}, Lemma 1). In order to present our main result, 
we introduce the following: 

\begin{definition}
Let $Y$ be an irreducible nodal curve, let $f: C \to Y$ be the normalization 
map and let $d$ be the gonality of $C$, so that there exists a morphism 
$\varphi: C \to \mathbb{P}^1$ of degree $d$. Assume that $\varphi$ is not 
composed with an involution, hence the locus $\Gamma := \{ (p,q) \in C \times C 
\setminus \Delta: \varphi(p)=\varphi(q) \}$ is irreducible (indeed, let 
$p \sim q$ if and only if $(p,q) \in \Gamma$, so that $\varphi$ factors 
as $C \to {C /\sim} \to \mathbb{P}^1$ and either the first map has degree $1$ 
and $\Gamma = \Delta$ or the second map has degree $1$ and $\Gamma$ is 
uniquely determined by $\varphi$). We say that a node $\chi$ on $Y$ is 
\emph{general $\varphi$-neutral} if $(p,q)$ is general in $\Gamma$ 
with $f^{-1}(\chi) = \{p, q \}$.   
\end{definition}

In particular, under the operation of making a general $\varphi$-neutral node 
$\varphi$ induces a pencil of degree $d$ on $Y$, while a general node increases 
by one the degree of the pencil induced by $\varphi$ on $Y$. 
We are able to prove that a curve obtained from a general $d$-gonal curve 
by making $n_1$ general nodes and $n_2$ general $\varphi$-neutral nodes exhibits 
the same Koszul cohomology vanishings as the general smooth $(d+n_1)$-gonal 
curve of the same genus. More precisely, the following holds: 

\begin{theorem}\label{main}
Fix integers $n_1 \ge 0$, $n_2 \ge 0$, $g \ge n := n_1+n_2$, 
and $d$ such that $2 \le d \le \lfloor (g-n+2)/2 \rfloor$. 
Let $C$ be a smooth $d$-gonal curve of genus $g-n$ 
such that $\dim G^1_{d+m}(C) \le m$ for all $m$ 
with $0 \le m \le g-n-2d+2$, let $\varphi: C \to \mathbb{P}^1$ 
be the degree $d$ morphism computing the gonality of $C$ and assume 
that $\varphi$ is not composed by an involution. Let $Y$ be a nodal 
curve of arithmetic genus $g$ with $C$ as its normalization, 
$n_1$ general nodes, and $n_2$ general $\varphi$-neutral nodes. 
Then $K_{g-d-n_1+1,1}(Y,\omega _Y)=0$ and $K_{r-d-n_1+1,1}(Y, L)=0$ for every 
line bundle $L$ on $Y$ with $h^0(Y, L)=r+1$ and $\deg(L) \ge 3g$.  
\end{theorem}

Notice that the above assumptions on $C$ hold for the 
general $d$-gonal curve (see \cite{a}, pp. 393--394). 
Recall also from \cite{bf1} the definition of the following locally 
closed algebraic subset of the moduli space $\overline{\mathcal{M}}_g$ 
of stable curves of genus $g$:
\begin{eqnarray*}
W(g,x,k,y) &:=& \{ X \in \overline{\mathcal{M}}_g: X \textrm{ is irreducible, }
\mathrm{card}(\mathrm{Sing}(X))=x, \\
& & \textrm{there exists a rank one torsion free sheaf $F$ on $X$}\\ 
& & \textrm{with } \deg(F)=k, h^0(X,F) \ge 2, \mathrm{card}(\mathrm{Sing}(F))=y \}.
\end{eqnarray*}
If $V(g,x,k,y)$ is the irreducible component of $W(g,x,k,y)$ whose general 
element has the general $(k-y)$-gonal curve as its normalization, then 
from Theorem \ref{main} we obtain the following:

\begin{corollary}
Let $Y$ be a general element of $V(g,x,k,y)$. 
If $2 \le k-y \le \lfloor (g-x+2)/2 \rfloor$, then
$K_{g-k+y+1,1}(Y,\omega _Y)=0$ and $K_{r-k+y+1,1}(Y, L)=0$ for every 
line bundle $L$ on $Y$ with $h^0(Y, L)=r+1$ and $\deg(L) \ge 3g$.
\end{corollary}

\section{The proofs}

\begin{lemma}\label{involution}
Let $C$ be an integral projective curve and $\varphi: C \to \mathbb{P}^1$, 
$\psi: C \to \mathbb{P}^1$ be morphisms. Assume $d:= \deg (\varphi) \ge 2$ 
and there is no morphism $i: \mathbb {P}^1 \to \mathbb {P}^1$ such that 
$\psi = i \circ \varphi$. Fix a general $P \in \mathbb {P}^1$. 
Then there are $p,q \in \varphi^{-1}(P)$ such that $\psi(p) \ne \psi(q)$.
\end{lemma}

\proof Assume that this is not true for a fixed $P \in \mathbb{P}^1$ 
such that $\mathrm{card}(\varphi^{-1}(P))=d$. Then $\psi(\varphi^{-1}(P))$ is a 
unique point, say $j(P)$. If the same holds for every sufficiently general $P$, 
then we get a rational map $j$ from $\mathbb{P}^1$ into itself such that 
$\psi = j \circ \varphi$ on a non-empty open subset of $C$. 
Since $\psi(\varphi^{-1}(P))$ is a unique point for a general $P\in \mathbb {P}^1$ 
and $\psi, \varphi$ are morphism, $\psi(\varphi^{-1}(P))$ is a unique point 
for every $P \in \mathbb {P}^1$. Hence $j$ may be uniquely extended to 
a set-theoretic map $i: \mathbb{P}^1 \to \mathbb{P}^1$. Since $\mathbb {P}^1$ 
is a smooth curve and $j$ is rational, $i$ is a morphism. 
Since $\psi, \varphi$ are morphisms and $C$ is separated, 
$\psi = i \circ \varphi $, contradiction.

\qed

\noindent \emph{Proof of Theorem \ref{main}.} 
Fix integers 
\begin{eqnarray*}
k &:=& g-d-n_1+1 \\
\nu &:=& g-2d-n_1+n_2+2. 
\end{eqnarray*}
and let $X$ be the stable curve obtained from $Y$ by identifying 
$\nu-n+1$ pairs of general points on $Y$. In particular, let $p, q$ 
be a pair of points on $Y$ identified to a node on $X$. If 
$K_{k,1}(Y, \omega_Y(p+q))=0$ then according to \cite{av}, 
Theorem 2.1, for every effective divisor $E$  of degree 
$e \ge 1$ we have $K_{k+e,1}(Y, \omega_Y(p+q+E))=0$. Thus 
if $L$ is any line bundle on $Y$ of degree $x \ge 3g$, 
then $h^0(Y,L-\omega_Y(p+q)) \ge 1$ and 
$K_{k+x-2g,1}(Y, L)=0$. On the other hand, 
by \cite{av}, proof of Lemma 2.3,
we have $K_{k,1}(Y, \omega_Y) \subseteq 
K_{k,1}(Y, \omega_Y(p+q)) \subseteq 
K_{k,1}(X, \omega_X)$, therefore in order 
to prove our statement we may assume 
$K_{k,1}(X, \omega_X) \ne 0$ and look for a contradiction.
By \cite{a}, Proposition 8, there exists a 
torsion-free sheaf $F$ on $X$ with $\deg(F)=k+1$ and $h^0(X,F) \ge 2$. 
Let $s$ with $0 \le s \le \nu+1$ be the number of nodes at which $F$ 
is not locally free. If $f: X' \to X$ is the partial 
normalization of $X$ at all such nodes, then $F=f_*(L)$, where 
$L = f^*(F) / \mathrm{Tors}(f^*(F))$ is a line bundle on $X'$ 
with $\deg L = k+1-s$ and $h^0(X',L)=h^0(X,F) \ge 2$. 
By taking the pull-back of $L$ on $C$, we obtain a 
$g^1_{k+1-s}$ not separating the $\nu+1-s$ pairs of points 
$(p_i,q_i)$ on $C$ glued to the nodes $\chi_i$ on $X$, 
$i=1, \ldots, \nu+1-s$. 

Assume first that the induced morphism $\psi: C \to \mathbb{P}^1$ 
of degree $k+1-s-b$ with $b \ge 0$ is equal to $\varphi$ composed 
with a morphism $\mathbb{P}^1 \to \mathbb{P}^1$. 
Since $\psi$ does not separate at least 
$\nu+1-s-n_2$ pairs of general points, if we let 
\begin{eqnarray*}
G^{\varphi}_{k+1-s}(C) &:=& \{ g^1_{k+1-s} \in G^1_{k+1-s}(C):
\textrm{the induced morphism is equal to}\\
& & \varphi \textrm{ composed with } i: \mathbb{P}^1 \to \mathbb{P}^1 \}, 
\end{eqnarray*}  
then we have
$$
\dim G^{\varphi}_{k+1-s}(C) \ge \nu-s-n_2+1.
$$
On the other hand, we have 
\begin{eqnarray*}
\dim G^{\varphi}_{k+1-s}(C) &=& \dim \{ i: \mathbb{P}^1 \to \mathbb{P}^1,
\deg(i) = \frac{k+1-s-b}{d} \}/ 
\mathrm{Aut}(\mathbb{P}^1) \\
&=& \frac{2}{d}(k+1-s-b)-2 \le \frac{2}{d}(k+1-s)-2 =\\
&=& \frac{2}{d}(\nu+d-n_2-s)-2=\frac{2}{d}(\nu-s-n_2)
< \nu-s-n_2+1
\end{eqnarray*}
for every $d \ge 2$, contradiction. 

Assume now that $\psi$ is not composed with a morphism 
$\mathbb{P}^1 \to \mathbb{P}^1$. Let us define inductively
\begin{eqnarray*}
G^{1,0}_{k+1-s}(C) &:=& \{ g^1_{k+1-s} \in G^1_{k+1-s}(C): 
\textrm{the induced morphism is not composed} \\
& & \textrm{ with } i: \mathbb{P}^1 \to \mathbb{P}^1 \}\\
G^{1,i}_{k+1-s}(C) &:=& \{ g^1_{k+1-s} \in G^{1,i-1}_{k+1-s}(C)
\textrm{ not separating } (p_i, q_i) \}.
\end{eqnarray*}  
We claim that $\dim G^{1,i}_{k+1-s}(C) \ge \nu+1-s-i$. 
Indeed, we have $G^{1,\nu+1-s}_{k+1-s}(C)\ne \emptyset$. 
Let $V$ be an irreducible component of $G^{1,i}_{k+1-s}(C)$,
assume by induction that $\dim V \ge \nu+1-s-i$ and let 
$W$ be the irreducible component of $G^{1,i-1}_{k+1-s}(C)$ 
containing $V$. If $\chi_i$ is a general node, then it is 
clear that $V \subsetneq W$, hence $\dim G^{1,i-1}_{k+1-s}(C)
\ge \dim W \ge \dim V + 1 \ge \nu+1-s-i+1$.
If instead $\chi_i$ is a general $\varphi$-neutral node and
$V=W$, then every linear series $g^1_{k+1-s}$ in $W$ induces a 
morphism $\psi: C \to \mathbb{P}^1$ such that $\psi(p_i)=\psi(q_i)$ 
for a general choice of $(p_i,q_i)$ with $\varphi(p_i)=
\varphi(q_i)$. By Lemma \ref{involution}, this is 
a contradiction, so the claim holds and in particular 
we have $\dim G^1_{k+1-s}(C) \ge \dim G^{1,0}_{k+1-s}(C) \ge
\nu+1-s$. 

In order to reach a contradiction, assume first 
$0 \le s \le g-n-2d+2$.
Hence we obtain   
$k+1-s-n_2 = d-2d+g-n+2-s$ with $0 \le -2d+g-n+2-s \le g-2d-n+2$
and our numerical hypotheses imply that  
$$
\dim G^1_{k+1-s-n_2}(C) \le g-2d-n+2-s.
$$
On the other hand, by \cite{fhl}, Theorem 1, we have 
$$
\dim G^1_{k+1-s}(C) \le \dim G^1_{k+1-s-n_2}(C)+2n_2,
$$
hence it follows that 
$$
\dim G^1_{k+1-s}(C) \le g-2d-n+2-s+2n_2 = \nu - s.
$$

Assume now $s > g-n-2d+2$. We claim that also in this case 
$$
\dim G^1_{k+1-s}(C) = \max_r \{ 2(r-1) + \dim W^r_{k+1-s}(C) \}
< \nu+1-s.
$$
Indeed, we have 
\begin{eqnarray*}
\dim W^r_{k+1-s}(C) &\le& \dim W^1_{k+1-s-(r-1)}(C) \le \\
&\le& \dim W^1_d(C) +2(k+1-s-(r-1)-d) \le \\
&\le& 1 +2(k+1-s-(r-1)-d) 
\end{eqnarray*}
where the second inequality is provided by \cite{fhl}, Theorem 1
(indeed, if $k+1-s-(r-1)<d$ then $W^1_{k+1-s-(r-1)}(C) = \emptyset$ 
since $C$ is $d$-gonal). 
Hence it follows that $\dim W^r_{k+1-s}(C) < \nu+1-s-2(r-1)$ 
for any $r$, as claimed. 

\qed

\providecommand{\bysame}{\leavevmode\hbox to3em{\hrulefill}\thinspace}

\end{document}